\documentclass[11pt]{amsart}
\usepackage[hidelinks]{hyperref}
\usepackage{doi,url}

\usepackage[utf8]{inputenc}
\usepackage{amsmath}
\usepackage{amssymb}
\usepackage[numbers]{natbib}

\usepackage{prettyref}
\usepackage{enumitem}
\usepackage{imakeidx}

\numberwithin{equation}{section}

\newtheorem{theorem}{Theorem}

\newtheorem{lemma}[theorem]{Lemma}
\newtheorem{fact}[theorem]{Fact}

\newcounter{introthmcounter}
\theoremstyle{definition}

\theoremstyle{remark}

\newcommand{\bE}{\mathbb E}

\newcommand{\bU}{\mathbb{U}}

\def\-{\text{-}}

\newcommand{\dcl}{\mathrm{dcl}}
\newcommand{\convex}{\mathrm{convex}}

\newcommand{\cO}{\mathcal{O}}

\newcommand{\RCF}{\mathrm{RCF}}

\DeclareMathOperator{\val}{\mathbf{v}}
\DeclareMathOperator{\res}{\mathbf{r}}

\DeclareMathOperator{\cof}{\mathrm{cof}}

\title{$T$-convexly valued o-minimal fields are definably spherically complete}
\author{Pietro Freni}


\begin{document}

\begin{abstract}
    I prove the statement in the title using results from \cite{freni2024t}. This shows that Question~1.1 in \cite{bradley-williams2023spherically} has negative answer for certain expansions of a valued field.
\end{abstract}

\maketitle
In \cite[Quest.~1.1]{bradley-williams2023spherically} the authors ask whether a \emph{definably spherically complete} expansion of a valued field necessarily has a spherically complete elementary extension. Here \emph{definably spherically complete} means that every definable nest of valuation balls has non-empty intersection. 

It is well known, by \cite{kuhlmann1997exponentiation}, that non-trivially $T$-convexly valued o-minimal fields defining exponentiation cannot be spherically complete. I will show that nevertheless they are still definably spherically complete.

Throughout $T$ will be any o-minimal theory expanding $\RCF$. If $\bE \models T$, and $S$ is a subset of an elementary extension of $\bE$, I will write $\bE \langle S \rangle:=\dcl_T(\bE\cup S)$. As usual $T_\convex$ will denote the common theory of all expansion $(\bE,\cO)$ of some $\bE\models T$ by a non-trivial $T$-convex valuation ring $\cO$ (see \cite{dries1995t}). Given $(\bE, \cO)\models T_\convex$, by $\bE$-definable I will mean definable with parameters in the o-minimal structure $\bE$ and by $(\bE, \cO)$-definable I will mean definable with parameters in the structure expanded by the valuation ring $\cO$. I will also write $\val(\bE,\cO)$ and $\res(\bE,\cO)$ for the imaginary sorts of the value group and the residue field respectively and consider them with their induced structure. The cofinality of a subset $X$ of an ordered set will be denoted by $\cof(X)$.

\begin{lemma}\label{lem:toname}
    If $(\bE, \cO)\models T_\convex$, then for every $(\bE, \cO)$-definable subset $X \subseteq \val(\bE, \cO)$, $\cof(X)\in \{\cof(\res(\bE, \cO)), \cof(\val(\bE, \cO)^{<0}),\cof(\bE),1\}$.
    \begin{proof}
        Let $b\in \bU \succ \bE \models T$ be such that $\cO <b<\bE^{>\cO}$, so that $\cO$ is externally definable using $b$. By \cite[(3.10)]{dries1995t}, every $(\bE, \cO)$-definable subset $X$ of $\bE^{>0}$ can be written as a boolean combination of intervals and preimages of $\cO$ by monotone $\bE$-definable functions. In particular such $X$ is the trace in $\bE$ of a set definable in $\bE\langle b\rangle$, and thus, by o-minimality a finite union of points and convex sets of the form $\bE \cap (a_0,a_1)$ where $a_0, a_1 \in \bE \langle b\rangle$. Now if $a_1 \in \bE$, then $\cof( (a_0, a_1) \cap \bE)= \cof(\bE)$. Otherwise $\cof( (a_0, a_1) \cap \bE) \in \{\cof(\bE^{<b}), \cof(\bE^{>b})\}$ and one easily reckons $\cof(\bE^{<b})= \cof(\res(\bE, \cO))$, $\cof(\bE^{>b})= \cof(\val(\bE, \cO)^{<0})$. Thus every $(\bE, \cO)$-definable subset of $\bE^{>0}$ has a maximum or one of these three cofinalities. Finally observe that if $X \subseteq \bE^{>0}$, then $\cof(\val(X)) \in \{1, \cof(\{1/x: x \in X\})\}$.
    \end{proof}
\end{lemma}

As in \cite[Def.~3.15]{freni2024t}, given $(\bE, \cO) \models T_\convex$ and an uncountable cardinal $\kappa$, I will call \emph{$\kappa$-bounded wim-constructible extension} an elementary extension of $(\bE, \cO)$ of the form $(\bE\langle x_i: i<\mu\rangle, \cO')$, where $\mu$ is an ordinal and $(x_i: i< \mu)$ is a $\mu$-tuple from some elementary extension of $(\bE,\cO)$ such that for every $j<\mu$, $x_j$ is a pseudolimit for some p.c.\ sequence of length $<\kappa$ in $\bE\langle x_i: i<j\rangle$ which has no pseudolimit in $\bE\langle x_i: i<j\rangle$.

\begin{fact}\label{fact:thmA}
    If $(\bE, \cO) \prec (\bE_*, \cO_*)$ is a $\kappa$-bounded wim-constructible extension then $\res(\bE, \cO)=\res(\bE_*, \cO_*)$, $\val(\bE, \cO)^{<0}$ is cofinal in $\val(\bE_*, \cO_*)^{<0}$ and $\bE$ is cofinal in $\bE_*$.
    \begin{proof}
        By \cite[Thm.~A]{freni2024t} if an element $x$ of an elementary extension $(\bE_*, \cO_*)$ of $(\bE, \cO)$ is a pseudolimit of a p.c.\ sequence in $(\bE, \cO)$ without pseudolimits therein, then $\res(\bE \langle x \rangle, \cO_* \cap \bE \langle x \rangle) = \res(\bE, \cO)$, $\bE$ is cofinal in $\bE\langle x \rangle$ and $\val(\bE, \cO)^{<0}$ is cofinal in $\val(\bE\langle x \rangle, \cO_* \cap \bE \langle x \rangle)^{>0}$. The statement then easily follows by transfinite induction.
    \end{proof}
\end{fact}

\begin{theorem}
    Every model of $T_\convex$ is definably spherically complete.
    \begin{proof}
        Let $(\bE, \cO)\models T_\convex$ and $\kappa$ be a cardinal greater than $\cof(\res(\bE, \cO))$, $\cof(\bE)$, and $\cof(\val(\bE, \cO)^{<0})$. Let $(\bE_*, \cO_*)\succ (\bE, \cO)$ be a maximal $\kappa$-bounded wim-constructible extension of $(\bE,\cO)$ and notice that in particular $(\bE_*, \cO_*)$ is $\kappa$-spherically complete. It suffices to observe that $(\bE_*, \cO_*)$ is definably spherically complete: by Fact~\ref{fact:thmA} we have $\cof(\res(\bE, \cO)) = \cof(\res(\bE_*, \cO_*))$, $\cof(\bE)=\cof(\bE_*)$ and $\cof(\val(\bE, \cO)^{<0})=\cof(\val(\bE_*, \cO_*)^{<0})$, thus by Lemma~\ref{lem:toname}, in $(\bE_*, \cO_*)$ every nested $(\bE_*, \cO_*)$-definable family of valuation balls has non-empty intersection.
    \end{proof}
\end{theorem}

\subsection*{Acknowledgments} I thank David Bradley-Williams for making me aware of the question and Vincenzo Mantova for a related discussion. This happened at the conference MAC30 held from the 7th to the 9th of Octber 2024 in Leeds in occasion of the 30th year of Dugald Macpherson at the university of Leeds and I thank the organizers of the conference.


\begin{thebibliography}{4}
\providecommand{\natexlab}[1]{#1}
\providecommand{\url}[1]{\texttt{#1}}
\expandafter\ifx\csname urlstyle\endcsname\relax
  \providecommand{\doi}[1]{doi: #1}\else
  \providecommand{\doi}{doi: \begingroup \urlstyle{rm}\Url}\fi

\bibitem[Bradley-Williams and Halupczok(2023)]{bradley-williams2023spherically}
D.~B. Bradley-Williams and I.~Halupczok.
\newblock Spherically complete models of {H}ensel minimal valued fields.
\newblock \emph{MLQ Math. Log. Q.}, 69\penalty0 (2):\penalty0 138--146, 2023.
\newblock ISSN 0942-5616,1521-3870.

\bibitem[Freni(2024)]{freni2024t}
P.~Freni.
\newblock T-convexity, weakly immediate types and {$T$}-{$\lambda$}-spherical
  completions of o-minimal structures, 2024.
\newblock URL \url{https://arxiv.org/abs/2404.07646}.

\bibitem[Kuhlmann et~al.(1997)Kuhlmann, Kuhlmann, and
  Shelah]{kuhlmann1997exponentiation}
F.-V. Kuhlmann, S.~Kuhlmann, and S.~Shelah.
\newblock Exponentiation in power series fields.
\newblock \emph{Proc. Amer. Math. Soc.}, 125\penalty0 (11):\penalty0
  3177--3183, 1997.
\newblock ISSN 0002-9939,1088-6826.
\newblock \doi{10.1090/S0002-9939-97-03964-6}.
\newblock URL \url{https://doi.org/10.1090/S0002-9939-97-03964-6}.

\bibitem[van~den Dries and Lewenberg(1995)]{dries1995t}
L.~van~den Dries and A.~H. Lewenberg.
\newblock {$T$}-convexity and tame extensions.
\newblock \emph{J. Symbolic Logic}, 60\penalty0 (1):\penalty0 74--102, 1995.
\newblock ISSN 0022-4812,1943-5886.
\newblock \doi{10.2307/2275510}.
\newblock URL \url{https://doi.org/10.2307/2275510}.

\end{thebibliography}
\end{document}